\documentclass[12pt]{amsart}
\usepackage{amsmath}
\usepackage{enumerate}
\usepackage{amssymb}
\usepackage{amsbsy}
\usepackage{amsfonts}
\usepackage[arrow]{xy}

\theoremstyle{plain}
\newtheorem{thm}{Theorem}[section]
\newtheorem{prop}[thm]{Proposition}
\newtheorem{lem}[thm]{Lemma}
\newtheorem{cor}[thm]{Corollary}
\theoremstyle{definition}
\newtheorem{dfn}[thm]{Definition}
\newtheorem{ex}[thm]{Example}
\newtheorem{rmk}[thm]{Remark}
\numberwithin{equation}{section}
\newcommand{\sm}{\left(\begin{smallmatrix}}
\newcommand{\esm}{\end{smallmatrix}\right)}

\newfont{\FieldFont}{msbm10 scaled\magstep1}

\newcommand{\Mj}{{\mathcal{M}^{(j)}}}

\newcommand{\pf}{\noindent\bf Proof }

\def\sb{\mathop{\frak b}}

\begin{document}

\title{ Period Relations,  Jacobi Forms and  Eichler Integral}

\author{YoungJu  Choie }

 \address{Department of Mathematics and PMI\\
 Pohang University of Science and Technology\\
 Pohang, 790--784, Korea}
 \email{yjc@postech.ac.kr}

\author{Subong Lim }

 \address{Department of Mathematics and PMI\\
 Pohang University of Science and Technology\\
 Pohang, 790--784, Korea}
 \email{subong@postech.ac.kr}

 \thanks{Keynote:  Eichler Integral, cusp forms,
 mock modular forms, mock Jacobi forms, harmonic  Maass forms, period }
 \thanks{This work
 was partially supported by KOSEF R01-2008-000-20448-0(2008)
 and KRF-2007-412-J02302}

\begin{abstract}
We study period relations of Jacobi forms.  It turns out that the
relations satisfied by  Mordell integral coming from  Lerch or
Appell sums are the special case of those. The existence of Jacobi
integral associated to given period function using generalized
Poincar\'e series is claimed.

\end{abstract}
 \maketitle

\today




\section{\bf{Introduction}}

It is shown, when  Zwegers studied Ramanujan Mock theta functions,
the    Mordell integral\cite{Mor},
$$h(\tau,z):=\int_{\mathbb{R}} \frac{e^{\pi i \tau x^2-2\pi x
z}}{\cosh{\pi x}} dx$$ satisfies the following relations:
$$-\frac{e^{\pi iz^2/\tau}}{\sqrt{-i\tau}}h(-\frac{1}{\tau},
\frac{z}{\tau})+h(\tau,z)=0,$$
$$ h(\tau,z)=e^{\frac{\pi i}{4}}h(\tau+1,z)+e^{-\frac{\pi i}{4}}
\frac{e^{\frac{\pi i z^2}{\tau+1}}}{\sqrt{\tau+1}}
h(\frac{\tau}{\tau+1}, \frac{z}{\tau+1}),$$
$$h(\tau,z)+h(\tau, z+1)=\frac{2}{\sqrt{-i\tau}} e^{\pi i
\frac{(z+\frac{1}{2})^2}{\tau}},$$ and
$$h(\tau,z)+e^{-2\pi i z-\pi i \tau} h(\tau, z+\tau)=2 e^{-\pi i z
-\pi i \frac{\tau}{4}}.$$
\\

It turns out that these are the part of period relations
associated to Jacobi forms, namely, any  period function
$P(\tau,z)$ of Jacobi integral of weight $k $ and index
 $m$ (with trivial multiplier system)
 satisfies
 $$P(\tau,z)+ \tau^{-k} e^{-2\pi i m{\frac{z^2}{\tau}}} P(-\frac{1}{\tau},
 \frac{z}{\tau}) =0,$$
 $$P(\tau,z)+ (\frac{\tau-1}{\tau})^{-k} e^{-2\pi i m
 \frac{z^2}{\tau}}  P(\frac{\tau-1}{\tau},
 \frac{ z}{\tau})+( \frac{ -1}{\tau-1})^{-k}
 e^{-2\pi i m \frac{z^2}{\tau-1}}   P( \frac{ -1}{\tau-1},
 \frac{z}{\tau-1})=0.$$
\\

  From the recent work by
Zwegers\cite{Ze}, Bringmann-Ono\cite{BO1} it turns out that the
mock theta functions, which were studied by Ramanujan in his
letter\cite{R}, are  holomorphic parts of weak Maass forms. Based
on the modular behavior of mock theta functions Zagier\cite{Za}
further defined a concept of mock modular forms. However mock
modular form can be considered as a special case of modular
integral with period.
\\

The concept of modular integral already was introduced by Eichler
and studied further by many researchers(see, for instance,
\cite{K, K1, K2}). It is well known that Eichler integral plays a
role to understand periods of modular forms, which are related to
the modular symbols and  special values of L-functions(see
\cite{KZ}). Note that  a connection between period and Maass wave
forms was explored by Lewis and Zagier\cite{LZ} and further
applications have been explored  by many researchers\cite{BO2, HM,
M, M1, Mu} since then.
\\

The purpose of this article is to study period relations by
introducing a concept of Jacobi integral.  In particular we
introduce a concept of mock Jacobi form, which was already
appeared in several places(see \cite{BZ,Ze}), that is a
holomorphic Jacobi integral with a "dual" (true) Jacobi
form\cite{CI-2}. It turns out that Lerch sums studied in \cite{Ze}
and Appell functions studied in \cite{STY} can be viewed as
typical examples of mock Jacobi forms.
\\

This paper is organized as follows. We introduce some useful
notations in section 2. In section 3, the concept of Jacobi
integral with period functions has been introduced and  a lifting
map from Jacobi integrals to Jacobi forms are studied. Examples
from the indefinite theta series, Appell function and Jacobi
Eisenstein series of weight $2$ are introduced. In section 4,
period relations, using the relations of Jacobi group, are derived
and it is also explained in terms of the parabolic cohomology in
the sense of Eichler cohomology\cite{K}.   A family of Jacobi
integral with theta decomposition was introduced.

In section5, using a generalized Jacobi Poincar\'e series the
existence of Jacobi integral, which may have poles, was
claimed.
Here we modify the idea by Knopp \cite{K}, that is,  to
introduce a generalized Poincar\'e series to study Eichler
cohomology. The detailed proof goes to Appendix in the final
section. In section 6, we study a "mock Jacobi form" and period
relations of a family of mock Jacobi forms. Section 7 gives a
conclusion of this paper.

\section{\bf{Definitions and Notations}}

Let us set up the following notations. Let  $\mathcal{H}$ be the
usual complex upper half plane and $\tau \in \mathcal{H}, z\in
\mathbb{C}^j, j\geq 1.$   $\Gamma:=\Gamma(1):=SL(2, \mathbb{Z}). $
The Jacobi group $\Gamma^J $ is defined as follows:

\begin{dfn} Let
$$ \Gamma^{J}:= \Gamma  \propto \mathbb{Z}^{2j}= \{ [M,(\lambda,\mu)]| M
\in \Gamma, \lambda,\mu \in \mathbb{Z}^{j} \}.$$
\end{dfn}

This set $\Gamma^J$ forms a group under a group law
$$[M_1,(\lambda_1,\mu_1)][M_2,(\lambda_2,\mu_2)]=[M_1M_2, (\lambda',\mu')
+(\lambda_2,\mu_2)],$$ where $\sm \lambda'\\ \mu' \esm = M_2^t\sm
\lambda_1\\ \mu_1\esm$ and is called the {\bf Jacobi group}.  Note
that the Jacobi group $\Gamma^{J}$ acts on $\mathcal{H}\times
\mathbb{C}^j$ as, for each  $\gamma=[\sm a&b\\c&d\esm, (\lambda,
\mu )]\in \Gamma^J, \lambda, \mu \in \mathbb{Z}^j,$
$$\gamma (\tau,z)=(\frac{a\tau+b}{c\tau+d}, \frac{z+\lambda\tau+\mu}{c\tau+d}).$$

Furthermore, for
\begin{math}
    \gamma = [\sm a&b\\c&d \esm , (\lambda, \mu)]
\end{math} $\in \Gamma^{ J},  k \in \frac{1}{2}\mathbb{Z} $
and $\mathcal{M}^{(j)}\in M_{j\times j}(\frac{1}{2}\mathbb{Z}),$
let
$$j_{k,\mathcal{M}^{(j)} }(\gamma,(\tau,z)) := (c\tau+d)^{-k}e^{2 \pi i
Tr(  \mathcal{M}^{(j)}(-z^t \frac{c}{ c\tau+d} z
 +\lambda^t \tau \lambda+2\lambda^t
z+\lambda^t\mu))}.$$


Let us define the usual slash operator on a function $f:
\mathcal{H}\times \mathbb{C}^j  \rightarrow \mathbb{C}$ :
 $$(f|_{\omega,  k,\mathcal{M}^{(j)}} \gamma) (\tau,z) :=
 \omega(\gamma)j_{k,\mathcal{M}^{(j)}}(\gamma,(\tau,z))f(\gamma(\tau,z)),\
\gamma \in \Gamma^J,$$ where $\omega(\gamma)$ is the multiplier
system of weight $k$ on  $\Gamma^J $ so that it  satisfies
$$\omega(\gamma_1 \gamma_2)
j_{k,\mathcal{M}^{(j)}}(\gamma_1\gamma_2,(\tau,z))=\omega(\gamma_1)\omega(\gamma_2)
j_{k,\mathcal{M}^{(j)}}(\gamma_1, \gamma_2(\tau,z))
j_{k,\mathcal{M}^{(j)}}(\gamma_2,(\tau,z)),$$ for all
$\gamma_1,\gamma_2\in \Gamma^J.$ Then one checks the following
consistency condition (see also \cite{EZ}, Section I.1 ):
$$
(f|_{\omega, k,\mathcal{M}^{(j)}} \gamma |_{\omega,
k,\mathcal{M}^{(j)}} \gamma')(\tau,z)=(f|_{\omega,
k,\mathcal{M}^{(j)}} \gamma \gamma')(\tau,z),\ \gamma, \gamma'\in
\Gamma^J.$$ Throughout this paper we let
$$f|_{\omega, k,\mathcal{M}^{(j)}}\gamma = f|_{\omega} \gamma$$ unless it is specified.
Also when $\omega$ is trivial, i.e. $\omega(\gamma)=1,$ for all
$\gamma\in \Gamma^J$ we denote it as
$$f|_{\omega,k,\mathcal{M}^{(j)}}\gamma=f|_{k,\mathcal{M}^{(j)}}\gamma = f|\gamma.$$
Throughout this paper we let $v:=  Im(\tau),\
y:=Im(z)=(Im(z_1),\cdots,Im(z_j)).$


\section{\bf{Jacobi Integral }}

\medskip

Let $\mathcal{M}^{(j)}$ be fixed and
$\mathcal{P}_{\mathcal{M}^{(j)}}$ be the space of functions $f$
holomorphic in $\mathcal{H}\times\mathbb{C}^j$ which satisfy the
growth condition
$$|f(\tau,z)| < K(|\tau|^\rho+v^{-\sigma})e^{2\pi Tr( \frac{
\mathcal{M}^{(j)} y^t y}{v})},\ $$ for some positive constants
$K,\rho$ and $\sigma$.

\begin{prop} The set $\mathcal{P}_{\mathcal{M}^{(j)}}$ has the following properties:
\begin{enumerate}
\item It is preserved under $|_{\omega, k,\mathcal{M}^{(j)}}$ for
any real $k$ and any $\gamma \in \Gamma^J$. \item It forms a
vector space over $\mathbb{C}.$
\end{enumerate}
\end{prop}

{\pf} For simplicity we may assume that $j=1$ and
$\mathcal{M}^{(j)}=m\in \frac{1}{2}\mathbb{Z}.$

\noindent (1) It is enough to check for $[S,(0,0)],[T,(0,0)]$
 and $[I,(\lambda,\mu)],$
 where $S = \sm 1&1\\0&1 \esm, T = \sm 0&-1\\1&0 \esm, I = \sm 1&0\\0&1
 \esm$ since $\Gamma(1)^J$ is generated by those elements (see Section \ref{per}):
\begin{enumerate}
\item[(a)] For $[S,(0,0)]$,
\begin{eqnarray*}
&&|(f|_{\omega,k,m}[S,(0,0)])(\tau,z)| =|f(\tau+1,z)|\\
&<& K(|\tau+1|^\rho+v^{-\sigma})e^{2\pi m\frac{y^2}v}\\
&<& K'(|\tau|^{\rho'}+v^{-\sigma'})e^{2\pi m\frac{y^2}v}
\end{eqnarray*}
for some positive constants $K',\rho'$ and $\sigma'$.

\item[(b)] For $[T,(0,0)]$,
\begin{eqnarray*}
&&|(f|_{\omega,k,m}[T,(0,0)])(\tau,z)| = |(\tau)^{-k}e^{2\pi
mi(\frac{-z^2}\tau)}f(-\frac1\tau,\frac z\tau)|\\
&<& |\tau|^{-k}K(|\tau|^{-\rho}+(\frac
v{|\tau|^2})^{-\sigma})|e^{2\pi
mi(\frac{-z^2}\tau)}|e^{2\pi m\frac{Im(z\bar{\tau})^2}{v|\tau|^2}}\\
&<& K'(|\tau|^{\rho'}+v^{-\sigma'})e^{2\pi m\frac{y^2}v}
\end{eqnarray*}
for some positive constants $K',\rho'$ and $\sigma'$.

\item[(c)] For $[I,(\lambda,\mu)]$,
\begin{eqnarray*}
&&|(f|_{\omega,k,m}[I,(\lambda,\mu)])(\tau,z)| = |e^{2\pi
mi(\lambda^2\tau+2\lambda z)}f(\tau,z+\lambda\tau+\mu)|\\
&<& K(|\tau|^\rho+v^{-\sigma})|e^{2\pi mi(\lambda^2\tau+2\lambda
z)}|e^{2\pi m\frac{(y+\lambda v)^2}v}\\
&<& K(|\tau|^\rho+v^{-\sigma})e^{2\pi m\frac{y^2}v}.
\end{eqnarray*}
\end{enumerate}
\medskip

\noindent (2) For any $f,g\in  \mathcal{P}_{m} $  note that
$f(\tau,z)e^{-2\pi m\frac{y^2}v}$ and $g(\tau,z)e^{-2\pi
m\frac{y^2}v}$ satisfy the following growth condition
\begin{equation}\label{growthcondition}
|f(\tau,z)e^{-2\pi m\frac{y^2}v} |  <
K_1(|\tau|^{\rho_1}+v^{-\sigma_1}),
\end{equation}
$$|g(\tau,z)e^{-2\pi m\frac{y^2}v} |  <
K_2(|\tau|^{\rho_2}+v^{-\sigma_2}),$$ for some positive constants
$K_i,\rho_i$ and $\sigma_i, i=1,2.$ We conclude that $f+g\in
\mathcal{P}_{m}$. {\qed}

\medskip

 More generally, let
$\mathcal{P}$ be the space of functions which are holomorphic
in $\mathcal{H}$ with the growth condition in
(\ref{growthcondition}).
 It is easy to verify that $\mathcal{P}$ is
preserved under $|_{\omega,k,0}$ for any  $\gamma\in \Gamma^J$ and
forms a ring.

\medskip


\begin{dfn}\label{integral}(Jacobi Integral)
\begin{enumerate}

\item A real analytic  periodic
 function $f: \mathcal{H} \times \mathbb{C}^j
\rightarrow \mathbb{C}$ is called a   {\bf{Jacobi Integral }} of
weight $k\in \frac{1}{2} \mathbb{Z}$ and index
$\mathcal{M}^{(j)}$ with multiplier system $\omega$ and a
holomorphic period functions $P_{\gamma}$ on $\Gamma^J$ if it
satisfies the following relations:

\begin{enumerate}

\item[(i)] For all $\gamma \in \Gamma^J$
\begin{equation}{\label{Int}}
(f|_{\omega, k, \mathcal{M}^{(j)}} \gamma)(\tau,z) = f(\tau, z)+
P_{\gamma}(\tau,z),
\end{equation}
where $P_{\gamma}$ is in $\mathcal{P}_{\mathcal{M}^{(j)}}.$

\item[(ii)] It satisfies a growth condition, when $v, y
\rightarrow \infty,$
$$|f(\tau,z)|v^{-\frac{k}{2}}
e^{2\pi Tr(\mathcal{M}^{(j)}
\frac{y^t y}{v})} \rightarrow 0.$$
\end{enumerate}

\item The space of Jacobi integrals forms a vector space over
$\mathbb{C}$ and we denote it as $J_{\omega, k,
\mathcal{M}^{(j)}}^{\int}(\Gamma^J).$ In particular when $j=1$ we
let $\mathcal{M}^{(j)}=m$ so the space is denoted by $J_{\omega,
k, m }^{\int}(\Gamma^J).$


\end{enumerate}

\end{dfn}

\begin{rmk}
\begin{enumerate}
\item The periodic condition on $f$ is equivalent to say that
$$f|_{\omega,k,\mathcal{M}^{(j)}}[I,(0, 1_n)]=f, f|_{\omega, k,\mathcal{M}^{(j)}}[S^\ell,(0,0)]=f,$$  or
$$P_{[I, 0,1_n ]}(\tau,z)=0, P_{[S^{\ell}, 0,0 ]}(\tau,z)=0,$$
for some $\ell\in \mathbb{Z}$ and for all $n = 1,2,\cdots,j$ where
$1_n \in \mathbb{Z}^j$ whose $n$th component is $1$ and all other
components are $0$.

\item The collection of holomorphic functions $\{ P_{\gamma}|
\gamma \in \Gamma^J \}$ occurring in (\ref{Int}) is called  the
system of { \bf{period functions}} of $f.$ The period functions
$P_{\gamma}$ satisfy the following consistency condition:
$$P_{\gamma_1 \gamma_2} = P_{\gamma_1 }|_{\omega,k,\mathcal{M}^{(j)}}\gamma_2 +
P_{\gamma_2},\ \text{for all}\ \gamma_1, \gamma_2\in \Gamma^J.$$

\item If $P_{\gamma}(\tau,z)=0,$ for all $\gamma\in \Gamma^J,$
then $f$  is a usual {\bf{Jacobi form}}, whose space will be
denoted by $J_{\omega, k, \mathcal{M}^{(j)}}(\Gamma^J).$
\end{enumerate}
\end{rmk}

\medskip

\subsection{Lifting from Jacobi integrals to Jacobi forms}
In this section we take $j=1$ and $\mathcal{M}^{(j)}=m$ for
simplicity. The result in this section can be extended to general
$j\geq 1$ without any technical difficulties.

 Let us define the following operator $\Psi$ on $J_{\omega,k,m}^{\int}(\Gamma^J)$
 as
$$\Psi(f )(\tau, z):=v^{k} e^{-4\pi  m \frac{y^2}{v}}
(\frac{\partial{f}}{\partial \overline{z}})(\tau,z).$$

 Then the following
holds:

\begin{prop}
 Let $G(-\overline{\tau}, \overline{z}):=\Psi(f(\tau,z))$ with
$f(\tau,z)\in  {J^{\int}_{k,m}}(\Gamma^J).$ Then $G(\tau,z)$ is in
$J_{1-k,-m}(\Gamma^J).$
\end{prop}

{\pf} (1) Let $\gamma=[\sm a&b\\c&d\esm,(0,0)]\in\Gamma^J.$ By the
definition of Jacobi integral, we have
$$(c\tau+d)^{-k}e^{2\pi im\frac{-cz^2}{c\tau+d}}f(\frac{a\tau+b}{c\tau+d},\frac
z{c\tau+d})=f(\tau,z)+P_\gamma(\tau,z).$$ Then since $f(\tau,z)$
is real analytic with respect to $z$, we see that
$$(c\tau+d)^{-k}e^{2\pi im\frac{-cz^2}{c\tau+d}}\frac{\partial
f}{\partial\bar{z}}(\frac{a\tau+b}{c\tau+d},\frac
z{c\tau+d})\frac1{c\bar{\tau}+d} = \frac{\partial
f}{\partial\bar{z}}(\tau,z).$$ From this it follows that
\begin{eqnarray*}
&&\Psi(f)(\frac{a\tau+b}{c\tau+d},\frac z{c\tau+d})\\
&=&\frac{v^k}{|c\tau+d|^{2k}}e^{-4\pi
m\frac{(Im(z(c\bar{\tau}+d)))^2}{v|c\tau+d|^2}}\frac{\partial
f}{\partial\bar{z}}(\frac{a\tau+b}{c\tau+d},\frac z{c\tau+d})\\
&=&\frac{v^k}{|c\tau+d|^{2k}}e^{-4\pi
m\frac{(Im(z(c\bar{\tau}+d)))^2}{v|c\tau+d|^2}}(c\tau+d)^k(c\bar{\tau}+d)e^{2\pi
im\frac{cz^2}{c\tau+d}}\frac{\partial
f}{\partial\bar{z}}(\tau,z)\\
&=&\frac{v^k}{(c\bar{\tau}+d)^{k-1}}e^{2\pi
im\frac{c\bar{z}^2}{c\bar{\tau}+d}}\frac{\partial
f}{\partial\bar{z}}(\tau,z)\\
&=& (c\bar{\tau}+d)^{1-k}e^{2\pi
im\frac{c\bar{z}^2}{c\bar{\tau}+d}}\Psi(f)(\tau,z).
\end{eqnarray*}
Since $G(-\bar{\tau},\bar{z}) = \Psi(f)( \tau,z)$, we have
\begin{eqnarray*}
&&G(\frac{a\tau+d}{c\tau+d},\frac z{c\tau+d})\\
&=&\Psi(f)(-\frac{a\bar{\tau}+b}{c\bar{\tau}+d},\frac{\bar{z}}{c\bar{\tau}+d})\\
&=&\Psi(f)(\frac{a(-\bar{\tau})-b}{-c(-\bar{\tau})+d},\frac{\bar{z}}{-c(-\bar{\tau})+d})\\
&=&(-c(-\tau)+d)^{1-k}e^{2\pi
im\frac{-cz^2}{-c(-\tau)+d}}\Psi(f)(-\bar{\tau},\bar{z})\\
&=&(c\tau+d)^{1-k}e^{2\pi i(-m)\frac{cz^2}{c\tau+d}}G(\tau,z).
\end{eqnarray*}
(2) Let $\gamma=[I,(\lambda,\mu)]\in\Gamma^J$. By the definition
of Jacobi integral, we have
$$e^{2\pi im(\lambda^2\tau+2\lambda z)}\frac{\partial
f}{\partial\bar{z}}(\tau,z+\lambda\tau+\mu) = \frac{\partial
f}{\partial\bar{z}}(\tau,z).$$ From this it follows that
\begin{eqnarray*}
&&\Psi(f)(\tau,z+\lambda\tau+\mu)\\
&=&v^ke^{-4\pi m\frac{(y+\lambda v)^2}v}\frac{\partial
f}{\partial\bar{z}}(\tau,z+\lambda\tau+\mu)\\
&=&v^ke^{-4\pi m\frac{(y+\lambda v)^2}v}e^{-2\pi
im(\lambda^2\tau+2\lambda z)}\frac{\partial
f}{\partial\bar{z}}(\tau,z)\\
&=&e^{-2\pi
im(\lambda^2\bar{\tau}+2\lambda\bar{z})}\Psi(f)(\tau,z).
\end{eqnarray*}
Since $G(\tau,z) = \Psi(f)(-\bar{\tau},\bar{z})$, we have
\begin{eqnarray*}
&&G(\tau,z+\lambda\tau+\mu)\\
&=&\Psi(f)(-\bar{\tau},\bar{z}+\lambda\bar{\tau}+\mu)\\
&=&\Psi(f)(-\bar{\tau},\bar{z}-\lambda(-\bar{\tau})+\mu)\\
&=&e^{-2\pi im(\lambda^2(-\tau)-2\lambda
z)}\Psi(f)(-\bar{\tau},\bar{z})\\
&=&e^{-2\pi i(-m)(\lambda^2\tau+2\lambda z)}G(\tau,z).
\end{eqnarray*} So the proof is completed. {\qed}

\medskip

\subsection{Examples}

We   give several examples of real analytic Jacobi forms, whose
holomorphic part or non holomorphic part can be regarded as Jacobi
integrals.

The first example is  from that in Zwegers\cite{Ze}:
\begin{ex}\label{Zeex}
For $z=x+iy \in \mathbb{C}$ and $\tau=u+iv \in \mathcal{H},$
consider  the series
$$R(\tau,z)=\sum_{\nu\in
\frac{1}{2}+\mathbb{Z}}\{
sgn(\nu)-E((\nu+\frac{y}{v})\sqrt{2v})\}(-1)^{\nu-\frac{1}{2}}
e^{-\pi i \nu^2\tau-2\pi i \nu z},$$
$$E(z)=sgn(z)(1-\beta(z^2)),\
\beta(x)=\int_{x}^{\infty}u^{-\frac{1}{2}} e^{-\pi u} \, du \,
(x\in \mathbb{R} \geq 0 ).$$

Take a multiplier system
\[\omega([I,(1,0)]) = \omega([I,(0,1)]) = -1,\]
\[\omega([S,(0,0)]) = e^{\frac{\pi i}4}, \omega([T,(0,0)]) = \frac{-1}{\sqrt{-i}}.\]

Then

\begin{enumerate} \item $(R|_{\omega, \frac{1}{2}, -\frac{1}{2}}[T,(0,0)])(\tau,
z)=R(\tau, z)+P_{[T,(0,0)]}(\tau,z),$ where
$$P_{[T,(0,0)]}(\tau,z)=\int_{\mathbb{R}}
\frac{e^{\pi i \tau x^2- 2 \pi z x}}{\cosh{\pi x}} \, dx.$$

 \item For $a\in (0,1), b\in \mathbb{R}$ and $\tau \in
\mathcal{H},$  let
  $$R_{a,b}(\tau)=-i\int_{-\overline{\tau}}^{i\infty}
\frac{g_{a,-b}(\tau)}{\sqrt{-i(z+\tau)}} dz,$$

where
$$g_{a,b}(\tau):=\sum_{\nu\in a+\mathbb{Z}} \nu e^{\pi i
\nu^2\tau+2\pi i \nu b}.$$

Then
$$R_{a,b}(\tau)=ie^{-\pi i(a-\frac{1}{2})^2\tau-2\pi
i(a-\frac{1}{2})b} R(\tau, (a-\frac{1}{2})\tau+b+\frac{1}{2}).$$

\item This is a real analytic Jacobi integral of weight
$\frac{1}{2}$ and index $-\frac{1}{2} $ with multiplier system
$\omega.$

\item Furthermore, $\Psi(R)(\tau,z)\ (=v^{\frac12}e^{2\pi\frac{y^2}v}(\frac{\partial R}{\partial\bar{z}})(\tau,z))
=\sqrt{2} \theta(-\overline{\tau}, \overline{z}),$ where
$\theta(\tau,z)$ is the well-known Jacobi Theta series defined as
$$\theta(\tau,z):=\sum_{\nu\in \frac{1}{2}+\mathbb{Z}} e^{\pi i
\nu^2\tau+2\pi i \nu(z+\frac{1}{2})}.$$

\end{enumerate}

\end{ex}

\begin{ex}
The following is a real analytic Jacobi Eisenstein
series\cite{C2}:
$$-\frac{1}{12} E_{2,1}^*(\tau,z)= \sum_{n,r\in \mathbb{Z}
\atop 4n\geq r^2} H(4n-r^2)q^n\xi^r+ \frac{2}{\sqrt{v}} \sum_{r,
f\in \mathbb{Z}} \beta(\pi f^2 v) q^{\frac{r^2-f^2}{4}}\xi^r $$

is a real analytic Jacobi form (Eisenstein series of weight $2$
and index $1$ on $\Gamma(1)^J$). Here, $H(n)$ denotes the Hurwitz
class number formula (see \cite{C2}).

\begin{enumerate}

\item Then $H^*(\tau,z):=\displaystyle\sum_{n,r\in \mathbb{Z}
\atop 4n\geq r^2}H(4n-r^2)q^n\xi^r$ is a (holomorphic) Jacobi
integral of weight $2$ and index $1$ with period function
\begin{equation}\label{period}
P_{[T,(0,0)]}(\tau,z)=\frac{\sqrt{2}}{8\pi i} \int_{0}^{i\infty}
\theta(t,0)(-i(t+\tau))^{-\frac{3}{2}} dt\cdot \theta(\tau, z).
\end{equation}

\item  $R(\tau,z) :=\frac{2}{\sqrt{v}} \displaystyle\sum_{r, f\in
\mathbb{Z}} \beta(\pi f^2 v) q^{\frac{r^2-f^2}{4}}\xi^r $ is a
 (real analytic) Jacobi integral of weight $2$ and index $1$
  with  the same period function $P_{[T,(0,0)]}(\tau,z)$  in (\ref{period}).

  \item It was shown\cite{HZ} that
$R(\tau,z)  =\frac{\sqrt{2}}{8\pi i}\int_{-\overline{\tau}}^{ i
\infty}
 \theta(t,0)
(-i(t+\tau))^{-\frac{3}{2}}
 dt \cdot \theta(\tau,z)$.

\end{enumerate}

\end{ex}

The Higher-Level Appell functions studied in \cite{STY} can be
regarded as Jacobi integral. For example it has the following
property:

\begin{ex}\label{appell}
Let
$$\mathcal{K}_1(q,x,y)=
\sum_{m\in\mathbb{Z}}\frac{q^{\frac{m^2}{2}}x^m}{1-xyq^m},
q=e^{2\pi i \tau},x=e^{2\pi i z}, y=e^{2\pi i w}, \tau \in
\mathcal{H}, z,w\in \mathbb{C}.$$

\noindent Further $G(\tau,z,w):= \frac{1}{\theta(\tau,z)}\cdot
\mathcal{K}_1(\tau,z,w)$ with $
\theta(\tau,z)=\displaystyle\sum_{\lambda \in
\mathbb{Z}}q^{\frac{\lambda^2}{2}}x^{\lambda}.$ Note that
$$\frac{e^{-\pi i \frac{z^2}{\tau}}}{\sqrt{-i\tau}}
\theta(-\frac{1}{\tau}, \frac{z}{\tau})=\theta(\tau,z).$$

Take a multiplier system
\[\omega([I,(1,0)]) = \omega([I,(0,1)]) = 1,\]
\[\omega([T,(0,0)]) = \sqrt{-i} \]
and let $\mathcal{M}^{(2)} = \sm 0&0\\0&-\frac12\esm$. Then it was
shown  that
$$(G|_{\omega, \frac{1}{2},\mathcal{M}^{(2)}}[T,(0,0)])(\tau,
z,w)=G(\tau, z,w)+P_{[T,(0,0)]}(\tau,z,w),$$ where
$$P_{[T,(0,0)]}(\tau,z,w)= -\frac{1}{2} e^{\pi i \frac{w^2}{\tau}}
\int_{\mathbb{R}} e^{-\pi x^2 } \frac{e^{-2\pi i x \frac{
w}{\sqrt{-i\tau}}}dx}{1-e^{-2\pi x \sqrt{-i\tau}}}.$$

So $G(\tau, z,w)$ is a Jacobi integral with weight $\frac{1}{2}$
and index $\frac{1}{2} \sm 1&0&\\0&-1\esm $ with a period function
$P_{[T,(0,0)]}(\tau,z,w).$

\end{ex}

\medskip

\section{\bf Period Relations}\label{per}

In this section we study more precise period relations of Jacobi
integral in $ J_{\omega, k, m}( \Gamma(1)^J).$

\subsection{Jacobi Group}\label{group}

Let us introduce the following notations:
$$G_0=[S,(0,0)],
G_1=[S,(1,0)], G_2=[T,(1,0)],$$
$$ G_3=[I,(1,0)], G_4=[I,(0,1)], I^J=[I, (0,0)],$$
$$ V=G_2^3G_1=[-TS,(1,-1)], R=G_2^3G_0=[-TS, (0, -1)].$$


We recall the following facts:

\begin{thm}
 \begin{enumerate}
\item $\Gamma(1)^J$ is generated by $G_1$ and $G_2.$

\item $\Gamma(1)^J$ is generated by $G_0$ and $G_2$.

\item $\Gamma(1)^J$ is generated by $G_2$ and $V.$ The generators
$G_2$ and $V$ satisfy the relations
$$G_2^4=V^3=I^{J},$$
$$VG_2^2=[I,(-1,-2)]G_2^2V=G_2^2[I,(1,2)]V=G_2^2V[I,(-2,-1)],$$
and these are the defining relations for $\Gamma(1)^J.$

\item $\Gamma(1)^J$ is generated by $G_2$ and $R.$

\end{enumerate}

\end{thm}
{\pf} See \cite{C-J}. {\qed}

\begin{cor}
  The generators $G_1$ and $G_2$
  of the group $\Gamma(1)^J$
satisfy the relations
$$G_2^4=(G_2^3G_1)^3=I^J,\
G_2^4=(G_2^3G_0)^3=I^J,\ G_2^4=R^3=I^J.$$

\end{cor}
{\pf} See \cite{C-J}. {\qed}

\begin{rmk} The following relations hold:
$$[I,(0,-1)]=G_0^{-1}G_2^{-2}G_0G_2^2,\
[I,(0,-1)]=G_0^{-1}G_2^{-2}G_0, $$
$$ [T,(0,0)]=[-I,(1,0)]G_2^3,\
[-I, (1,0)]=G_0^{-1}G_2^{-2}G_0.$$
\end{rmk}

\begin{dfn}{(\bf{parabolic element})}
We call  any element of the form  $[\sm 1& \ell\\0 & 1\esm,$
$(0,r)]$, $\ell, r\in \mathbb{Z},$ a {\em {"parabolic element"}}.

\end{dfn}

\medskip





\subsection{Period functions}

Classically, there are two relations period polynomial $p(\tau)$
associated with elliptic cusp forms  of weight $2-k$ should
satisfy (see \cite{KZ}), namely, $$
p(\tau)+\tau^{-k}p(-\frac{1}{\tau})=0$$ and
$$p(\tau)+\tau^{-k}p(\frac{\tau-1}{\tau})+(\tau-1)^{-k}p(\frac{-1}{\tau-1})=0.$$
\medskip

In this section we study the relations in which the period
function $P_{[T,0,0]}(\tau,z)$ associated with Jacobi integrals should
satisfy. In particular when $z=0$ we recover those period
relations from elliptic modular forms.

For simplicity we consider the case when  $j=1$ and
$\mathcal{M}^{(j)}=m.$

\begin{prop}\label{modular}
  The transformation formulas  of Jacobi integral  on $\Gamma(1)^J$ in
$(\ref{period})$  can be reduced to the following two relations:

\begin{enumerate}
\item $ f|_{\omega,k,m} [S, (0,0)]=f$,

 \item $ f|_{\omega,k,m}[T,(1,0)]=f+P_{[T,(1,0)]},   $ with
 $P_{[T,(1,0)]}\in \mathcal{P}_{m}.$

\end{enumerate}
\end{prop}

{\pf} Since $\Gamma(1)^J$ is generated by $G_0$ and $G_2$ the
result follows. {\qed}

\begin{thm}\label{relation1}({\bf{Period Relations}})
Take a multiplier system $\omega$ with  $\omega(-I)=1.$ If a
Jacobi integral $f$ is even and periodic with respect to $z$,
i.e., $P_{[-I,(0,0)]}(\tau,z) = P_{[I,(0,1)]}(\tau,z) = 0,$ then
the period functions $ P_{[T,(0,0)]}(\tau,z)$ and
$h(\tau,z):=P_{[I,(1,0)]}$ $(\tau,z)$ satisfy the following
properties:

\begin{enumerate}

\item\label{tr} $ P_{[T, (1,0)]}=P_{[T,(0,0)]} $.

\item  $P_{[T,(0,0)]} + P_{[T,(0,0)]}|_{\omega,k,m} [T,(0,0)]=0$.

\item $P_{[T,(0,0)]} +  P_{[T,(0,0)]}|_{\omega,k,m} [ST,(0,0)]+
P|_{\omega,k,m}[(ST)^2,(0,0)]=0 $.

\item $h |_{\omega,k,m}[T,(0,0)]= -P_{[T,(0,0)]} + P_{[T,(0,0)]}
|_{\omega,k,m}[I,(0,-1)] $.

\item $h  = P_{[T,(0,0)]} - P_{[T,(0,0)]}|_{\omega,k,m}[-I,(1,0)]
$.

\item $h  + h |_{\omega,k,m}[-I, (1,0)]=0$.

\end{enumerate}
\end{thm}
\medskip

{\pf}Note that
$P_{[\gamma,(\lambda,\mu)]} = P_{[-\gamma,(\lambda,\mu)]}$
since $f$ is even with respect to $z.$

\begin{enumerate}

\item
$f|_{\omega,k,m}[I,(0,-1)]|_{\omega,k,m}[T,(0,0)]=f+P_{[T,(0,0)]}=f+P_{[T,(1,0)]}$.

\item It follows from   $[T,(0,0)]^2=[-I,(0,0)]$.

\item It follows from  $[ST,(0,0)]^3 = [-I,(0,0)]$.

\item From $[I,(1,0)][T,(0,0)]=[T,(0,-1)]=[T,(0,0)][I,(0,-1)]$ it
follows that
$$h|_{\omega,k,m}[T,(0,0)] + P_{[T,(0,0)]}=P_{[I,(0,-1)]}.$$

\item It follows from $[T, (1,0)]=[T,(0,0)][I,(1,0)]$ and
(\ref{tr}).

\item By the definition of $h$
\begin{eqnarray*}
h|_{\omega,k,m}[-I,(1,0)] &=& P_{[I,(1,0)]}|_{\omega,k,m}[-I,(1,0)]\\
&=& P_{[-I,(0,0)]} - P_{[-I,(1,0)]} = -h. {\qed}
\end{eqnarray*}
\end{enumerate}

\medskip

\begin{ex} The following examples are in Example \ref{Zeex}: We set
$$P(\tau,z) := P_{[T,(0,0)]}(\tau,z)=\int_{\mathbb{R}}
\frac{e^{\pi i \tau x^2- 2 \pi z x}}{\cosh{\pi x}} \, dx$$
and
$$h(\tau,z) := P_{[I,(1,0)]}(\tau,z) = 2e^{-\pi iz-\pi i\tau/4}.$$

Then we have
\begin{enumerate}
\item $P(\tau,z)-(P|_{\omega}[I,(0,1)])(\tau,
z)$\\
= $P(\tau,z) + P(\tau,z+1) =
\frac{2}{\sqrt{-i\tau}}e^{\pi(z+\frac12)^2/\tau} =
(h|_{\omega}[T,(0,0)])(\tau,z+1).$

\item $P(\tau,z)-(P|_{\omega}[I,(1,0)])(\tau, z)\\
=P(\tau,z) + e^{-2\pi iz-\pi i\tau}P(\tau,z+\tau) =2e^{-\frac{\pi i \tau}{4}-\pi i
z} = h(\tau,z)$.

\item $P(\tau,z)+(P|_{\omega}[T,(0,0)])(\tau,z) $\\
$=P(\tau,z)-\frac{1}{\sqrt{-i\tau}}
e^{2\pi i \frac{1}{2} \frac{ z^2}{\tau}}
 P(\tau,z)=0.$

\item $(P|_{\omega}[T,(0,0)])(\tau,z)+
(P|_{\omega}[S,(0,0)])(\tau,z)
 + (P|_{\omega}[STS,(0,0)])(\tau,z)$\\
    $= -P(\tau,z) + e^{\frac{\pi i}4}P(\tau+1,z) + e^{-\frac{\pi i}4}
    \frac{e^{\pi iz^2/(\tau+1)}}{\sqrt{\tau+1}}P(\frac\tau{\tau+1},\frac z{\tau+1}) = 0$.
\end{enumerate}
\end{ex}

\begin{ex} The following examples are in Example \ref{appell}: We set
$$P(\tau,z,w) := P_{[T,((0,0),(0,0))]}(\tau,z,w)=e^{\pi i \frac{w^2}{\tau}}
\Phi(\tau,w) $$and
$$h(\tau,z,w) := P_{[I,((0,0),(1,0))]}(\tau,z,w) = e^{-2\pi iw-\pi i\tau}.$$
Here,   $$\Phi(\tau,w):= -\frac{1}{2} \int_{\mathbb{R}} e^{-\pi
x^2 } \frac{e^{-2\pi i x \frac{ w}{\sqrt{-i\tau}}}dx}{1-e^{-2\pi x
\sqrt{-i\tau}}}.$$  Then we have
\begin{enumerate}
\item $P(\tau,z,w)-(P|_{\omega}[I,((0,0),(0,-1))])(\tau,
z,w)$\\
= $e^{\pi\frac{w^2}\tau}\Phi(\tau,z,w) - e^{\pi i\frac{(w-1)^2}\tau}\Phi(\tau,w-1)\\
 = -\frac{i}{\sqrt{-i\tau}}e^{\pi i\frac{(w-1)^2}{\tau}}=-h|_w[T,((0,0),(0,0))].$

\item $P(\tau,z,w)-(P|_{\omega}[I,((0,0),(1,0))])(\tau, z,w)\\
=e^{\pi i\frac{w^2}\tau}\Phi(\tau,z,w)-e^{\pi i\frac{w^2}\tau}\Phi(\tau,z,w+\tau)
 = e^{-2\pi iw-\pi i\tau} = h(\tau,z,w)$.

\item $P(\tau,z,w)+(P|_{\omega}[T,((0,0),(0,0))])(\tau,z,w) $\\
$=e^{\pi i\frac{w^2}\tau}\Phi(\tau,z,w) + \frac{\sqrt{-i}}{\sqrt{\tau}}
\Phi(-\frac1\tau,\frac z\tau,\frac w\tau) = -1$.
\end{enumerate}
\end{ex}

\begin{prop}
 Proposition \ref{relation1}-(2) implies that
 Proposition \ref{relation1}-(5)
\end{prop}

{\pf} From the relation $[ST,(1,0)] [ST, (0,0)]
[ST,(0,0)]=[-I,(0,-1)]$ we derive
\begin{eqnarray*}
h|_{\omega,k,m}[-I,(0,0)]&=& -P|_{\omega,k,m}[ST,(0,0)]
-P|_{\omega,k,m}[(ST)^2,(0,0)]\\
&&-P|_{\omega,k,m}[I,(1,0)][-I,(0,0)]=0. {\qed}
\end{eqnarray*}

\medskip

Using the above relation we study a family of Jacobi integral
which has a theta decomposition

\subsection{ A  Jacobi integral with Theta decomposition}

Consider a holomorphic Jacobi integral
$\phi_{k,m}:\mathcal{H}\times \mathbb{C}^{2 } \rightarrow
\mathbb{C}$ such that

\begin{enumerate}\label{condition}
\item[(A)]$(\phi_{k,m}|_{\omega, k, \mathcal{M}}[S, ((0,0 ),(0,0
))])(\tau, z, w)=\phi_{k,m}(\tau,z,w)$.

\item[(B)]$(\phi_{k,m}|_{\omega, k, \mathcal{M}}[I, ((1, 0),
(1,0))])(\tau,z,w)  =\phi_{k,m}(\tau,z,w).$

\item[(C)] $(\phi_{k,m}|_{\omega, k, \mathcal{M}}[I, ((0, 0),
(0,1))])(\tau,z,w)  =\phi_{k,m}(\tau,z,w)$.
\end{enumerate}
Then the following holds:

\begin{prop}
\begin{enumerate}

\item $\phi_{k,m}(\tau,z,w)=
\displaystyle\sum_{\ell\pmod{2m}}h_{\ell}(\tau,w)\theta_{m,\ell}(\tau,z),$

\noindent where $h_{\ell}(\tau,w)=e^{-\frac{\pi i \ell^2 \tau}{2m}
} \int_{p}^{p+1}\phi(\tau,z,w)e^{-2\pi i \ell z} dz, p \in
\mathbb{C}$ and
$$\theta_{m,\ell}(\tau,z):=
\sum_{ r \in {\mathbb{Z}} \atop r\equiv \ell (mod 2m) }
q^{\frac{r^2}{4m}} \xi^{r}, q=e^{2\pi i \tau}, \xi=e^{2\pi i z}.$$

\item $(\phi_{ k,m}|_{\omega,k,\mathcal{M}}
[T,((0,0),(0,0))])(\tau,z,w)$ $= \phi_{
k,m}(\tau,z,w)+\displaystyle\sum_{\ell=0}^{2m-1}P_{\ell}(\tau,
w)\theta_{m,\ell}(\tau, z),$ where
$P_{\ell}(\tau,w)=h_{\ell}(\tau,w)-(h_{\ell}|_{\omega
k-\frac{1}{2},\mathcal{M}}[T, (0,0)])(\tau,w).$

\item $(\phi_{ k, m}|_{\omega, k,
\mathcal{M}}[I,((0,1),(0,0))])(\tau,z,w)$ $= \phi_{k,m}(\tau,z,w)
+\displaystyle\sum_{\ell=0}^{2m-1}(P_{\ell}(\tau,w) $
$$-(P_{\ell}|_{
\omega, k-\frac{1}{2},\mathcal{M} }[-I, ((1,0),
(0,0))])(\tau,w))\theta_{m,\ell}(\tau,z).$$

\end{enumerate}
\end{prop}

{\pf} $(1)$ The condition $(B)$ implies that
$$e^{2\pi i Tr(\mathcal{M}  ( \sm 1\\0\esm \tau (1,0) + 2 \sm 1\\0\esm  (z,w) )}
\phi_{k,m}(\tau,z+\tau,w)=\phi_{k,m}(\tau,z,w) $$ so that it has
the following theta series expansion, $$\phi_{k,m}(\tau,z,w)=
\displaystyle\sum_{\ell\pmod{2m}}h_{\ell}(\tau,w)\theta_{m,\ell}(\tau,z)
$$ (see detailed proof in \cite{EZ} or \cite{Ze}).
\\

\noindent $(2)$ This follows from transformation formula of Theta
Series (see \cite{EZ})
 $$ \theta_{m,\mu}(-\frac{1}{\tau},
\frac{z}{\tau}) =\sqrt{\frac{\tau}{2mi}} e^{2\pi i m
\frac{z^2}{\tau}} \sum_{\nu\pmod{2m}}e^{-2\pi i\frac{ \mu\nu}{2m}}
\theta_{m,\nu}(\tau,z).$$
\\

\noindent $(3)$ This follows from the period relation given in
Proposition \ref{relation1}: if $$ (\phi_{ k, m}|_{\omega, k,
\mathcal{M}}[I,((0,1),(0,0))])(\tau,z,w) =
\phi_{k,m}(\tau,z,w)+P_{[I,((0,1),(0,0))]}(\tau,z,w)$$ then
$$
P_{[I,((0,1),(0,0))]}(\tau,z,w)=P_{[T,((0,0),(0,0))]}(\tau,z,w)$$
$$-(P_{[T,((0,0),(0,0))]}|_{\omega,
k-\frac{1}{2},\mathcal{M} }[-I, ((1,0), (0,0))])(\tau,z,w).$$

{\qed}


\subsection{Cohomology}


We call any collection of functions $\{\varphi_\gamma|\ \gamma\in\Gamma^J\}$
in $\mathcal{P}_{m}$ which satisfies
\begin{equation}\label{cocycle}
\varphi_{\gamma_1\gamma_2} = \varphi_{\gamma_1}|_{\omega,k,m}\gamma_2 +
\varphi_{\gamma_2},\ \text{for}\ \gamma_1,\gamma_2\in\Gamma^J \end{equation} a
{\em{cocycle}} of weight $k$ and index $m$ on $\Gamma^J.$  A
{\em{coboundary}} of weight $k$ and index $m$ on $\Gamma^J$ is a
{\em{cocycle }} $\{\varphi_\gamma|\ \gamma\in\Gamma^J\}$ such that
$\varphi_\gamma = \varphi|_{\omega,k,m}\gamma - \varphi$ for all
$\gamma\in\Gamma^J,$  with $\varphi$ a fixed function in
$\mathcal{P}_{m}$. The parabolic cocycles on $\Gamma^J$ are the
cocycles $\{\varphi_\gamma|\ \gamma\in\Gamma^J\}$ which satisfy the
following additional condition: for each parabolic element $Q_j,
j=0,...,\ell,$ there exist $\varphi_j,
 \in \mathcal{P}_{m}, j=0,..,\ell,$ such that
\begin{eqnarray} \label{paraboliccocycle}
\varphi_{Q_j} &=& \varphi_j|_{\omega, k,m}Q_j -\varphi_j.
\end{eqnarray}

\begin{dfn} \begin{enumerate}
\item The cohomology group $H^1_{\omega, k,m}(\Gamma,
 \mathcal{P}_{m})$ is defined to be the vector space of cocycles
modulo coboundaries.

\item Let $\tilde{H}^1_{ \omega, k,m}(\Gamma,  \mathcal{P}_{m})$
be the subgroup of $H^1_{\omega, k,m}(\Gamma,  \mathcal{P}_{m})$
defined as {\em{the space of parabolic cocycles }} modulo
coboundries and we call $\tilde{H}^1_{ \omega, k,m}(\Gamma,
\mathcal{P}_{m})$ a parabolic cohomology group.

\end{enumerate}
\end{dfn}

\begin{rmk}
\begin{enumerate}
\item This is an analogous definition of the Eichler (parabolic)
cohomoloy group $ \tilde{H}^1_{  k, v  }(\Gamma, P_k),$ where
$P_k$ is the vector space of polynomials of degree $\leq k$ (see
\cite{K}).

\item For each $\phi\in J_{k,m}(\Gamma)$ there are at least  two
ways to attach the elements in $\tilde{H}^1_{ \omega, k ,
m}(\Gamma,  \mathcal{P}_{m}).$ One is via Eichler integral with
$\mathbb{Q}$-division points and the other is via Eichler Integral
and theta decomposition.

\end{enumerate}
\end{rmk}

Now take $\Gamma^J=\Gamma(1)^J$ and consider the space of the
following period functions:
\begin{eqnarray}\label{rr}
\mathcal{P}er_{\omega, k,m} &:=& \{ P\in  \mathcal{P}_{m} \, |
\,
P+P|_{\omega,k,m}[T,(0,0)]\\
\nonumber && =
P+P|_{\omega,k,m}[ST,(0,0)]+P|_{\omega,k,m}[(ST)^2,(0,0)]=0\}.
\end{eqnarray}

Then the following is true:

\begin{prop}
$\mathcal{P}er_{\omega, k,m}$ is a generating set of all parabolic
cocycles of $\Gamma(1)^J.$
\end{prop}

{\pf} First of all, one may regard $P \in  \mathcal{P}er_{\omega,
k,m}$ as the element  $P= P_{[T,(0,0)]} $ in the set of parabolic
cocycles on $\Gamma(1)^J$ from Proposition \ref{relation1}. On the
other hand,  the set  $\{[S,(0,0)], [T,(0,0)], [I,(0,1)]\}$
generates $\Gamma(1)^J $  and, again from Proposition
 \ref{relation1}, one has
$$P_{[I,(0,1)]}=P_{[T,(0,0)]} -
P_{[T,(0,0)]}|_{\omega,k,m}[-I,(1,0)].$$ So it is clear that
$P_{[T,(0,0)]}$   generates every parabolic cocycle in
$\Gamma(1)^J.$ {\qed}

\medskip


\section{\bf{  Jacobi Poincar\'e series and Existence of Jacobi Integral}}

\subsection{Jacobi Poincar\'e series}

A generalized Poincar\'e series was studied in \cite{K} to show
the isomorphism between the parabolic cohomology group  and space
of elliptic modular  cusp forms of the arbitrary weight.

Here we also introduce a generalized (Jacobi) Poincar\'e series to show the
existence of Jacobi integral, which may has poles,  associated to
given period functions $P_\gamma \in  \mathcal{P}_{m}$ on
$\Gamma(1)^J$.

\begin{dfn} Suppose $\{\varphi_\gamma \, | \, \gamma\in \Gamma(1)^J\}$ is a
parabolic cocycle of weight $k$ and index $m$ which satisfies the
additional condition that $\varphi_{[S,(0,0)]} = \varphi_{[I,
(0,1)]} = 0$. Suppose $k$ is a positive even integer and $\omega$ is
a multiplier system of weight $k$. Then the generalized Poincar\'e
series $\Phi(\{\varphi_\gamma\},k,m,\omega;\tau,z) = \Phi(\tau,z)$ is
defined by
\begin{equation} \label{series}
\Phi(\tau,z) := \displaystyle\sum_{\gamma \in \Gamma(1)_{\infty}^J
\backslash \Gamma(1)^J}
(\varphi_{\gamma}|_{\omega,k,m}\gamma)(\tau,z),
\end{equation}
where $\Gamma(1)_\infty^J = \{[\pm\sm 1&n\\0&1\esm,(0,\mu)]|\
n,\mu\in\mathbb{Z}\}$.
\end{dfn}

 Note that the assumption $\varphi_{[S,(0,0)]} =
\varphi_{[I, (0,1)]} = 0$ has been made to insure that the
individual terms of the series are independent of the choice of
$\gamma$ of coset representatives.

\begin{thm}\label{poin}
 For   sufficiently large $k$ $(>g)$ and $m>0$ the generalized
Poincar\'e series $\Phi(\{\varphi_\gamma\},k,m,\omega;\tau,z)$
converges absolutely where $g=max(2e+5,4),$ where $e$ is defined
in Lemma \ref{secondlemma} in Appendix.

\end{thm}
{\pf}
 The proof of absolute convergence of the series defining
$\Phi(\tau,z)$ is based upon a series of lemmas. Those are
essentially Lemma $4$, Lemma $5$ and Lemma $6$ in \cite{K}. We
give the detailed proof in the Appendix. {\qed}

\medskip

\subsection{  Existence of Jacobi Integral }

\begin{thm} \label{existence}
Let $r$ be any real number, $m>0$ and $\omega$ a multiplier system
of weight $r$ and index $m$. Suppose that $\{\varphi_ \gamma |\
\gamma\in\Gamma(1)^J\}$ is a parabolic cocycle of weight $k$ and index
$m$ in $\mathcal{P}_{m}$ such that $\varphi_{[S,(0,0)]} =
\varphi_{[I,(0,1)]} = 0$. Then there is a meromorphic function
$f:\mathcal{H}\times\mathbb{C}\rightarrow\mathbb{C}$ such that
$$(f|_{\omega, r,m}\gamma)(\tau,z) - f(\tau,z) = \varphi_{\gamma}(\tau,z)\
\text{for all}\ \gamma\in\Gamma(1)^J.$$
\end{thm}

{\pf} We take a generalized Poincar\'e series
$\Phi(\{\varphi_\gamma\},k,m,\omega;\tau,z) = \Phi(\tau,z)$ for
sufficiently large $k$ such that $k>g,$ $g=max(2e+5,4).$  For $\gamma
\in\Gamma(1)^J$, we see that
\begin{equation}\label{ftnaleqn}
\Phi|_{\omega, r,m} \gamma =
\bar{\omega}(\gamma)j_{k,m}(\gamma,(\tau,z))^{-1}\Phi(\tau,z) -
\bar{\omega}(\gamma)j_{k,m}(\gamma,(\tau,z))^{-1}g(\tau,z)\varphi_\gamma(\tau,z),
\end{equation}
 where
$g(\tau,z)$ is the Eisenstein series
\begin{equation} \label{Eisen}
g(\tau,z) := \displaystyle\sum_{\gamma \in \Gamma(1)_{\infty}
\backslash \Gamma(1)^J} \omega(\gamma)j_{k,m}(\gamma,(\tau,z)).
\end{equation}
The functional equation (\ref{ftnaleqn}) is a straightforward
consequence of the absolute convergence of (\ref{series}), the
consistency condition for the cocycle $\{\varphi_\gamma\}$ and the
consistency condition for the multiplier system $\omega$. For
$k\geq4$ the series in (\ref{Eisen}) converges absolutely and it
follows that
$$g(\gamma(\tau,z)) = \bar{\omega}(\gamma)j_{k,m}(\gamma,(\tau,z))^{-1}g(\tau,z)$$
for $\gamma\in\Gamma^J$. Thus, putting $F(\tau,z) :=
-\frac{\Phi(\tau,z)}{g(\tau,z)}$ and applying (\ref{ftnaleqn}), we
find that
\begin{eqnarray*}
(F|_{\omega, r,m} \gamma)(\tau,z) &=& -\frac{(\Phi|_{\omega, r,m} \gamma)(\tau,z)}{
g(\gamma(\tau,z))}\\
&=& -\frac{\Phi(\tau,z)}{g(\tau,z)} + \varphi_\gamma(\tau,z) =
F(\tau,z) + \varphi_\gamma(\tau,z)
\end{eqnarray*}
so that $F$ is a solution of the functional equation. {\qed}

\begin{rmk}
In \cite{K} the generalized Poincar\'e series has been studied to
show the isomorphism between the Eichler (parabolic) cohomology
group $\tilde{H}^1_{k,v}(\Gamma, P_k)$ and the space of cusp forms
of weight $2-k$ on $\Gamma,$ where the Petersson's result has been
used to guarantee that one can construct a modular form which has
the assigned  poles and zeros in $\mathcal{H}$ (see \cite{K} for
details).  However it is not known yet if the analogous result of
Petersson can be extended to the Jacobi form case to show the
constructed function $f$ in Theorem \ref{existence} is holomorphic
in $\mathcal{H}\times \mathbb{C}.$
\end{rmk}

\medskip







\section{\bf{  Mock Jacobi forms   }}\label{mock}


In this section we introduce  a mock Jacobi form, which has a
corresponding dual Jacobi form. Further study see \cite{CI-2}.

\subsection{Mock modular  form}

The concept of Mock modular form, which was motivated from
Ramanujan Mock Theta function,  was first introduced by Zagier in
\cite{Za}: A   function $H:\mathcal{H} \rightarrow \mathbb{C}$ is
called a {\bf{mock modular form }} if
\begin{enumerate}
\item It is holomorphic in $\mathcal{H}$ with only possible poles
at the cusps (so that it contains the weakly holomorphic modular
forms).

\item There is a rational number $\lambda$ such that $H(q),
q=e^{2\pi i \tau},$ must be multiplied by $q^{\lambda}$ in order
to have any kind of modularity properties, and a "shadow"
$g=\mathcal{S}[h]$ which is an ordinary modular form of weight
$2-k $ such that the holomorphic function
$h(\tau)=q^{\lambda}H(q)$ becomes a non-holomorphic modular form
of weight $k$ when we complete it by adding a correction term $g^*
(\tau)$ associated to $g(\tau).$

\item This "shadow" depends $\mathbb{R}$-linearly on $h$ and
vanishes if and only if $h$ is a modular form, so that we have an
exact sequence over $\mathbb{R}:$
$$0\rightarrow  {\mathbf{M}}_k^{!} \rightarrow \mathbb{M}_k
\rightarrow^{\mathcal{S}} M_{2-k}$$ Here,   $\mathbf{M}_k^{!},
\mathbb{M}_k$ and $M_{2-k}$ are the space of weakly holomorphic
modular forms, the space of mock modular forms and the space of
modular forms, respectively.
\end{enumerate}

\begin{rmk}
A mock modular form defined here is more restricted than that in
\cite{Za} since we take a rational invariant $\lambda=0$ (see
\cite{Za} for more detailed information).

\end{rmk}

The various examples were discussed by Zagier\cite{Za}. Here is
one more example, which was already computed in \cite{P}.

\begin{ex}\begin{enumerate}
\item Assume $k \in 2\mathbb{Z} $  and let $\tau=u+iv \in
\mathcal{H}.$  Consider
$$G_k(\tau|s)=\frac{1}{2} \sum_{m,n\in
\mathbb{Z}}^{'}
\frac{1}{(m\tau+n)^k}
(\frac{v}{|m\tau+n|^2}
)^s$$
which converges for $Re(s)>1-\frac{k}{2}. $

Then the Fourier expansions of $G_k(\tau|s)$ was derived in
\cite{P}:
$$G_k(\tau|s)=\zeta(k+2s)v^s+(-1)^{\frac{k}{2}} \pi
\frac{2^{2-k-2s} \Gamma(k-1+2s) }{ \Gamma(k+s)\Gamma(s) }
\zeta(k-1+2s) v^{1-k-s} $$
$$ +\frac{(-1)^{\frac{k}{2}}
(2\pi)^{k+2s}}{ \Gamma(k+s) \Gamma(s) }
 v^s \sum_{n\geq 1}
\sigma_{k-1+2s}(n)\{ \sigma(4\pi n v, k+s,s)e^{2\pi i n \tau}
$$ $$+\sigma(4\pi n v, s, k+s)e^{-2\pi i n \overline{\tau}} \}$$
 where
$\sigma_{\omega}(n)=\displaystyle\sum_{d|n\atop d>0} d^{\omega}$
and $\sigma(\eta, \alpha,
\beta)=\int_{0}^{\infty}(u+1)^{\alpha-1}u^{\beta-1}e^{-\eta u} \,
du$.

\medskip

 Let $G_k^*(\tau|s)=\pi^{-s}\Gamma(s)G_k(\tau|s)$.
Then
$$G_k^*(\tau|0)=H_{k}^*(\tau)+R_k^*(\tau),$$
where
$$H^*_k(\tau)=\frac{(-k)!}{ (2\pi
i)^{-k}} \zeta(1-k)-\frac{(1-k)!}{ (2\pi i)^{1-k}}
\int_{\tau}^{i\infty}
 [G_{2-k }(w|0)-\zeta(2-k)](w-\tau)^{-k} dw $$
$$ =\frac{(-k)!}{(2\pi i)^{-k}}
(\zeta(1-k)+\sum_{n\geq 1}
\sigma_{k-1}(n) e^{2\pi i n})$$ and
\begin{eqnarray*}
&&R_k^*(\tau)=(-k)!\zeta(2-k) (\frac{v}{ \pi})^{1-k}\\
&-& \frac{(1-k)!}{(2\pi i)^{1-k}}
\int_{-\overline{\tau}}^{i\infty}
[G_{2-k}(w|0)-\zeta(2-k)](w+z)^{-k} dw.
\end{eqnarray*}

\noindent So $H^*_k(\tau)$ is a mock modular form with Shadow
$\mathcal{S}[H_k^*(\tau)]=G_{2-k}(\tau|0).$
\end{enumerate}
\end{ex}

\bigskip

\subsection{Mock Jacobi forms}
Let us recall the following heat operator introduced in \cite{CK}:
Take a matrix $\Mj \in M_{j\times j}(\mathbb{R}).$ The heat
operator $L_{\Mj}$ is defined by
$$L_{\Mj}=8\pi i |\Mj|\frac{\partial}{\partial
\tau}-(\frac{\partial}{\partial z})^t \widetilde{\Mj} (
\frac{\partial}{\partial z}),$$ where
$$(\frac{\partial}{\partial z})=(\frac{\partial}{\partial
z_{\ell}}), z=(z_{\ell})_{1\leq \ell \leq j} \in \mathbb{C}^j $$
and $|\Mj|$ is the determinant of $\Mj, \widetilde{\Mj}=(
\widetilde{\Mj}_{mn}), \widetilde{\Mj}_{mn}$ is the cofactor of the
$(m,n)$th entry of $\Mj$ for $j\geq 2,$ and $\widetilde{\Mj}=1$ when
$j=1.$

\begin{dfn}\label{weak} A mock Jacobi form
 $\phi:\mathcal{H} \times \mathbb{C}^j\rightarrow \mathbb{C} $  is
 a meromorphic Jacobi integral in
 $J^{\int}_{\omega,-k+\frac{j}{2},\mathcal{M}^{(j)}}(\Gamma^J), k\in \mathbb{Z}_{\geq 0},$ such that
 $L_{\Mj}^{k+1}(\phi)$ is
  a  nontrivial (meromorphic)
Jacobi form of weight $k+\frac{j}{2}+2$ and index
$\mathcal{M}^{(j)}$ with multiplier system.  The Jacobi form
$L_{\Mj}^{k+1}(\phi)$ is called a "dual" of $\phi.$  In other
words, we say that a meromorphic  Jacobi integral which has a
"dual" Jacobi form is   a mock Jacob form.

\end{dfn}



\medskip

The following was introduced by Zwegers\cite{Ze}:
\subsection{Lerch Sum}

Consider the Lerch sum,
 $$\mu(\tau,z,w):=\frac{e^{\pi i w}}{\theta(\tau,z)}\sum_{n\in
\mathbb{Z}}\frac{(-1)^n e^{\pi i(n^2+n)\tau+2\pi i n z}}{1-e^{2\pi
i n \tau+2\pi i w}}, $$ which was originally studied by Lerch and
whose elliptic and modular transformation properties were derived
by Zwegers\cite{Ze} to connect with Mock theta function. Here
$$\theta(\tau,z):=\sum_{\nu \in \frac{1}{2}+\mathbb{Z}} e^{\pi i
\nu^2\tau+2\pi i \nu(z+\frac{1}{2})}.$$

Using the transformation properties of $\mu(\tau,z,w)$ and
$\theta(\tau,z)$ the followings were derived in \cite{Ze} and
\cite{BZ}:

\begin{ex}
\begin{enumerate}
\item Let $f(\tau,z) := e^{\pi iz-\pi
i\tau/4}\mu(\tau,z,\frac12\tau+\frac12)$ which is a Jacobi
integral of  a weight $\frac{1}{2}$ and an index $m=-\frac12.$ It
satisfies the following transformation properties:
\begin{enumerate}
\item
 \[\frac{-i}{\sqrt{-i\tau}}e^{\pi\frac{z^2}\tau}f(-\frac1\tau,\frac z\tau)
 = f(\tau,z) - \frac1{2i}e^{\pi iz-\pi i\tau/4}h(\tau,z-\frac12\tau+\frac12).\]
So, this implies that the period
 \[P(\tau,z) := P_{[T,(0,0)]}(\tau,z) = - \frac1{2i}e^{\pi iz-\pi i\tau/4}h(\tau,z-\frac12\tau+\frac12).\]

\item $L_m(f(\tau,z))$ is a "dual" of $f$, that is, a (nontrivial)
(meromorphic) Jacobi form of weight $\frac{5}{2}$ and index
$-\frac{1}{2}.$ Here,
\[L_m:=4\pi i\frac{\partial}{\partial\tau}+\frac{\partial^2}{\partial
z^2}\] is the corresponding heat operator.
\end{enumerate}

 \item  More
generally,let
  $f_{a,b}(\tau,z) := e^{2\pi iaz-\pi
  ia^2\tau}\mu(\tau,z,a\tau+b), $ for any $a,b\in \mathbb{R},$
  which is a Jacobi integral of
  weight $\frac{1}{2} $ and index $-\frac{1}{2} $ with its dual $L_m(f_{a,b})
  $ which is a Jacobi form of weight $\frac{5}{2}$ and index
  $-\frac{1}{2}.$

\item In fact that the dual of $f_{a,b}$  was computed explicitly
in \cite{BZ}:
 \[(4\pi i\frac{\partial}{\partial\tau}+\frac{\partial^2}
 {\partial z^2})(f_{a,b}(\tau,z))=e^{2\pi iaz-\pi ia^2\tau}
 \frac{16\pi^2\eta(\tau)^6}{\theta(\tau,a\tau+b)\theta(\tau,z)^3}\]
 \[\times
 \{\alpha_1(\tau)\theta_0(2\tau,2z+a\tau+b)-\alpha_0(\tau)\theta_1(2\tau,2z+a\tau+b)\}.\]
Here,
 \begin{eqnarray*}
 \theta_0(\tau,z) &:=& \sum_{n\in\mathbb{Z}}e^{\pi in^2\tau+2\pi inz},\\
 \theta_1(\tau,z) &:=& \sum_{n\in\frac12+\mathbb{Z}}e^{\pi in^2\tau+2\pi inz},\\
 \alpha_0(\tau) = \alpha_0^{a,b} &:=& \sum_{n\in\mathbb{Z}}(n+a/2)e^{2\pi in^2\tau+2\pi in(a\tau+b)},\\
 \alpha_1(\tau) = \alpha_1^{a,b} &:=& \sum_{n\in\frac12+\mathbb{Z}}(n+a/2)e^{2\pi in^2\tau+2\pi in(a\tau+b)}.
 \end{eqnarray*}

\end{enumerate}
\end{ex}

\begin{rmk} Further a family examples via Eichler integrals are constructed in \cite{CI-2}
\end{rmk}





\section{\bf{Conclusion}}
In this paper we study the period relations associated with Jacobi
integral.  This explains the relations from the Modell integral
associated to  Lerch sums\cite{Ze} and from the functional
relations associated  higher Appell functions\cite{STY}. On the
other hand, modular symbols can be studied purely algebraically
using period relations\cite{M} and recently modular symbols are
extended to the complex weight forms associated to Maass wave
forms.  We are intending to develop higher modular symbols,
Jacobi-modular symbols, using multi-variable period relations as
well as the actions of Hecke operators on them\cite{C3}.

\section{\bf{Appendix}}

Here we begin to prove Theorem \ref{poin}:

\begin{lem} \label{firstlemma}
For real numbers $c,d$ and $\tau=u+iv$, we have
$$(\frac{v^2}{1+4|\tau|^2})(c^2+d^2)\leq|c\tau+d|^2\leq2(|\tau|^2+v^{-2})(c^2+d^2).$$
\end{lem}

 If $A\in\Gamma(1)$ consider a factorization of $A$, $A =
C_1\cdots C_q$ where each $C_i$ is $T$ or a power of $S$. Eichler
showed that for any $A\in\Gamma(1)$ the factorization can be carried
out so that
$$q\leq m_1log\mu(A)+m_2,$$
where $m_1,m_2>0$ are independent of $A$ and
$$\mu(A) = a^2 + b^2 + c^2 + d^2\ \text{if}\ A = \sm a&b\\c&d
\esm.$$
 We assume that the cocycle $\{\varphi_\gamma\}$ in $\mathcal{P}_{m}$
satisfies
\begin{eqnarray} \label{generator}
|\varphi_{[T,(0,0)]}(\tau,z)| &<& K(|\tau|^\rho+v^{-\sigma})e^{2\pi m\frac{y^2}v},\\
\nonumber |\varphi_i(\tau,z)| &<&
K(|\tau|^\rho+v^{-\sigma})e^{2\pi m\frac{y^2}v},\ \text{for}\
0\leq i \leq l.
\end{eqnarray}

 Here $\varphi_i(\tau,z)$ is defined by (\ref{paraboliccocycle})
and $K,\rho,\sigma$ are positive constants. Assume also
$2\sigma>k, \rho>-k$.

\begin{lem} \label{secondlemma}If $\{\varphi_\gamma\}$ is a parabolic cocycle then there
exists $K^*>0$ such that
$$|\varphi_{[C_h,(0,0)]}|_{\omega,k,m}[C_{h+1},(0,0)]\cdots [C_q,(0,0)](\tau,z)|\leq
K^*\mu(A)^e\{|\tau|^{6e+2k}+v^{-6e-2k}\}e^{2\pi m\frac{y^2}v},$$ for
$1\leq h\leq q$. Here $e = max(\frac\rho2,\sigma-\frac k2)$ and
$A=C_1\cdots C_q$ is a factorization of $A\in\Gamma(1)$.
\end{lem}

{\pf} Consider first the case when $C_h$ is $T$. Let
$\gamma=C_{h+1}\cdots C_q = \sm a&b\\ c&d \esm$. Then,
by (\ref{generator}),
\begin{eqnarray*}
&&|\varphi_{[C_h,(0,0)]|[\gamma,(0,0)]}(\tau,z)| =
|c\tau+d|^{-k}|e^{2m\pi i(\frac{-c
z^2}{c\tau+d})}||\varphi_{[C_h,(0,0)]}(\gamma\tau,\frac
z{c\tau+d})|\\
&<&
|c\tau+d|^{-k}K(|\gamma\tau|^\rho+v^{-\sigma}|
c\tau+d|^{2\sigma})e^{2m\pi\frac{y^2}v}\\
&=&
(K|a\tau+b|^\rho|c\tau+d|^{-k-\rho}+K|
c\tau+d|^{2\sigma-k}v^{-\sigma})e^{2m\pi\frac{y^2}v}.
\end{eqnarray*}
By Lemma \ref{firstlemma},
\begin{eqnarray*}
|a\tau+b|^\rho &\leq&
2^{\frac\rho2}(|\tau|^2+v^{-2})^{\frac\rho2}(a^2+b^2)^{\frac\rho2},\\
|c\tau+d|^{2\sigma-k} &\leq& 2^{\sigma-\frac
k2}(|\tau|^2+v^{-2})^{\sigma-\frac
k2}(a^2+b^2)^{\sigma-\frac k2},\\
|c\tau+d|^{-k-\rho} &\leq& 2^{\sigma-\frac
k2}(|\tau|^2+v^{-2})^{\frac{-k-\rho}2}(a^2+b^2)^{\frac{-k-\rho}2}.
\end{eqnarray*}
Hence we have
\begin{eqnarray*}
|\varphi_{[C_h,(0,0)]}|[\gamma,(0,0)](\tau,z)| &<&e^{2\pi
m\frac{y^2}v}\{K2^{\frac\rho2}(|\tau|^2+v^{-2})^{\frac\rho2}
(a^2+b^2)^{\frac\rho2}\\
&\ &\times(\frac{1+4|\tau|^2}{v^2})^{\frac{\rho+k}2}(c^2
+d^2)^{\frac{-k-\rho}2}\\
&\ &+K2^{\sigma-\frac k2}(|\tau|^2+v^{-2})^{\sigma-\frac
k2}(c^2+d^2)^{\sigma-\frac k2}v^{-\sigma}\}.
\end{eqnarray*}
Since the nonzero $c$, $\sm *&*\\ c&* \esm\in\Gamma(1)$,
with $\Gamma(1)$ discrete, have a positive lower bound, it follows
that $c^2+d^2$ has a positive lower bound; hence
\begin{eqnarray*}
|\varphi_{[C_h,(0,0)]}|[\gamma,(0,0)](\tau,z)| &<&e^{2\pi
m\frac{y^2}v}\{K_1(a^2+b^2)^{\frac\rho2}(|\tau|^2
+v^{-2})^{\frac\rho2}(\frac{1+4|\tau|^2}{v^2})^{\frac{\rho+k}2}\\
&\ &+K_1'(c^2+d^2)^{\sigma-\frac
k2}(|\tau|^2+v^{-2})^{\sigma-\frac k2}v^{-\sigma}\}.
\end{eqnarray*}
Note that
$$a^2+b^2+c^2+d^2 = \mu(\gamma) \leq K_2\mu(A),$$
so that
\begin{eqnarray*}
|\varphi_{[C_h,(0,0)]}|[\gamma,(0,0)](\tau,z)| &\leq&e^{2\pi
m\frac{y^2}v}
\{K_3\mu(A)^{\frac\rho2}(|\tau|^2+v^{-2})^{\frac\rho2}v^{-k-\rho}
(1+4|\tau|^2)^{\frac{\rho+k}2}\\
&\ &+K_3'\mu(A)^{\sigma-\frac
k2}v^{-\sigma}(|\tau|^2+v^{-2})^{\sigma-\frac k2}\}.
\end{eqnarray*}
Letting $e=max(\frac\rho2,\sigma-\frac k2)$, we have
\begin{eqnarray*}
|\varphi_{[C_h,(0,0)]}|[\gamma,(0,0)](\tau,z)| &\leq& e^{2\pi
m\frac{y^2}v}\{K_4\mu(A)^e(|\tau|^2+v^{-2})\\
&\ &\times(v^{-k-\rho}(1+4|\tau|^2)^{\rho+\frac
k2}+v^{-\sigma})\}\\
&\leq& e^{2\pi
m\frac{y^2}v}\{K_4\mu(A)^e(|\tau|^2+v^{-2})^e\\
&\
&\times(\frac12v^{-2k-2\rho}+\frac12(1+4|\tau|^2)^{\rho+k}+v^{-\sigma})\}.
\end{eqnarray*}
Now $\sigma\leq e+\frac k2$ and $\rho+k\leq 2e+k$, so that
\begin{eqnarray*}
|\varphi_{[C_h,(0,0)]}|[\gamma,(0,0)](\tau,z)| &\leq& e^{2\pi
m\frac{y^2}v}\{K_5\mu(A)^e(|\tau|^2+v^{-2})^e(|\tau|^{4e+2k}+v^{-4e-2k})\}\\
&\leq& e^{2\pi
m\frac{y^2}v}\{K_6\mu(A)^e(|\tau|^{6e+2k}+v^{-6e-2k})\}.
\end{eqnarray*}
We now deal with the case in which $C_h = S^m$ for some
$m\in\mathbb{Z}$. Then
$$\varphi_{[S,(0,0)]} = \varphi_0|[S,(0,0)]-\varphi_0,$$
and therefore
$$\varphi_{[C_h,(0,0)]} = \varphi_0|[C_h,(0,0)]-\varphi_0.$$
From this it follows that
$$\varphi_{[C_h,(0,0)]}|[C_{h+1}\cdots C_q,(0,0)] =
\varphi_0|[C_h\cdots C_q,(0,0)]-\varphi_0|[C_{h+1}\cdots
C_q,(0,0)].$$ The previous argument applies to each of the two
terms on the righthand side to yield
$$|\varphi_{[C_h,(0,0)]}|[C_{h+1}\cdots C_q,(0,0)](\tau,z)|\leq
e^{2\pi m\frac{y^2}v}\{K_7\mu(A)^e(|\tau|^{6e+2k}+v^{-6e-2k})\}.$$
The proof is completed. {\qed}

\medskip

 Now we will use the Ford fundamental region $\mathcal{R}$. It is defined as follows:
$$\mathcal{R} =\{\tau\in\mathcal{H}|\ u<\frac\lambda2\ \text{and}\
|c\tau+d|>1\ \text{for all}\ \gamma=\sm
*&*\\c&d\esm\in\Gamma-\Gamma_\infty\}.$$ Then there exists $y_0>0$
with $iy_0\in\mathcal{R}$. Now determine $M$ by the
condition that $[A,(\lambda,0)]\in M$ if
$-\frac\lambda2\leq Re\{A(iy_0)\}<\frac\lambda2$.

\begin{lem} \label{thirdlemma}
If $[A,(\lambda,0)]=[\sm a&b\\c&d\esm,(\lambda,0)]\in M$,
chosen as indicated above, then
$$\mu(A)\leq K'(c^2+d^2),$$
for $K'>0$, independent of $A$.
\end{lem}

\begin{lem}\label{elliptic} For $\lambda\in\mathbb{Z}$
$$|\varphi_{[I,(\lambda,0)]}(\tau,z)| < |\lambda|K(|\tau|^\rho+v^{-\sigma})e^{2\pi
m\frac{y^2}v}.$$
\end{lem}

{\pf} First note that $\varphi_{[I,(0,0)]} = 0$ and
$$\varphi_{[I,(0,0)]} =
\varphi_{[I,(-1,0)]}|[I,(1,0)]+\varphi_{[I,(1,0)]}.$$ So
$|\varphi_{I,(-1,0)]}(\tau,z)| = |\varphi_{I,(1,0)]}(\tau,z)|.$
And for $\lambda>0$ we see that
$$|\varphi_{[I,(\lambda,0)]}(\tau,z)| =
|\varphi_{[I,(\lambda-1,0)]}|[I,(1,0)](\tau,z) +
\varphi_{[I,(1,0)]}(\tau,z)|.$$ So by induction on $\lambda$ we
get the result in the case of $\lambda>0$. And we can prove the
result when $\lambda<0$ by the same way. {\qed}

\begin{lem} \label{convergence}
\begin{enumerate}
\item The series
\begin{equation}\label{series1}
\sum_{c,d\in\mathbb{Z}\atop
(c,d)=1}\sum_{\lambda\in\mathbb{Z}}(c\tau+d)^{-k}e^{2\pi
m(\lambda^2\frac{a\tau+b}{c\tau+d}+2\lambda\frac
z{c\tau+d}-\frac{cz^2}{c\tau+d})}
\end{equation}
converges absolutely if $k>3$.

\item The series
\begin{equation} \label{series2}
\sum_{c,d\in\mathbb{Z}\atop
(c,d)=1}\sum_{\lambda\in\mathbb{Z}}(c\tau+d)^{-k}\lambda e^{2\pi
m(\lambda^2\frac{a\tau+b}{c\tau+d}+2\lambda\frac
z{c\tau+d}-\frac{cz^2}{c\tau+d})}
\end{equation}
converges absolutely if $k>4$.
\end{enumerate}
\end{lem}

{\pf} (1) We consider the series
$$\displaystyle\sum_{x\in\mathbb{Z}}e^{-\alpha(x+\beta)^2}.$$
Note that the following estimate holds:
\begin{equation*}
\sum_{x\in\mathbb{Z}}e^{-\alpha(x+\beta)^2}\leq
1+2\sum_{x\in\mathbb{N}}e^{-\alpha x^2}.
\end{equation*}
This is clear for $\beta\in\mathbb{Z}$. And if
$\beta\not\in\mathbb{Z}$, it follows from
\begin{displaymath}
-(x+\beta)^2 \leq \left\{
        \begin{array}{ll}
            -(x+[\beta]+1)^2 & \textrm{for $x\leq-[\beta]-1$, }\\
            -(x+[\beta])^2 & \textrm{for $x\geq-[\beta]$,}
        \end{array} \right.
    \end{displaymath}
where $[\beta]$ is the Gauss bracket. Since $e^{-\alpha x^2}>0$ and
this is a decreasing function we get the following estimate
$$\displaystyle\sum_{x\in\mathbb{N}}e^{-\alpha x^2}\leq\int^\infty_0e^{-\alpha x^2}dx.$$
And if we use
$$\int^\infty_0e^{-\alpha x^2}dx = \frac12\sqrt{\frac\pi \alpha}$$
we get
\begin{equation} \label{estimate}
\sum_{x\in\mathbb{Z}}e^{-\alpha(x+\beta)^2}\leq 1+ \sqrt{\frac\pi
\alpha}.
\end{equation}
 Now we estimate the sum
\begin{eqnarray*}
&\ &\sum_{\lambda\in\mathbb{Z}}|e^{2\pi
m(\lambda^2\frac{a\tau+b}{c\tau+d}+2\lambda\frac z{c\tau+d}-\frac{cz^2}{c\tau+d})}|\\
&=&|e^{2\pi
m(\frac{-cz^2}{c\tau+d})}|\sum_{\lambda\in\mathbb{Z}}e^{-2\pi
m(\lambda^2\frac
v{|c\tau+d|^2}+2\lambda\frac{Im(z(c\bar{\tau}+d))}{|c\tau+d|^2})}\\
&=&e^{2\pi m\frac{y^2}v}\sum_{\lambda\in\mathbb{Z}}e^{-2\pi m\frac
v{|c\tau+d|^2}(\lambda+\frac{Im(z(c\bar{\tau}+d))}v)^2}.
\end{eqnarray*}
If we use (\ref{estimate}) then we see that the series
(\ref{series1}) converges absolutely if $k>3$.

(2) Note that if $Re(\alpha)<0$ then we have
$$\int^\infty_{-\infty}xe^{\alpha x^2+\beta x}dx =
\frac{\sqrt{\pi}}{\sqrt{-\alpha}}(-\frac
\beta{2\alpha})e^{-\frac{\beta^2}{4\alpha}}$$ and
$$|Im(z(c\bar{\tau}+d))| \leq |z||c\tau+d|.$$
Then by the same argument we see that the series (\ref{series2})
converges absolutely if $k>4$. {\qed}

\medskip

Note that the series (\ref{series1}) can be written as
$$\displaystyle\sum_{[A,(\lambda,0)]\in\Gamma(1)_{\infty}^J \backslash \Gamma(1)^J}1|_{k,m}[A,(\lambda,0)]$$
and the series (\ref{series2}) can be written as
$$\displaystyle\sum_{[A,(\lambda,0)]\in\Gamma(1)_{\infty}^J \backslash \Gamma(1)^J}
\lambda(1|_{k,m}[A,(\lambda,0)]).$$

{\pf of Theorem \ref{poin}} Suppose $[A,(\lambda,0)]\in M$.
As before write $A=C_1\cdots C_q$. Then we find that
\begin{eqnarray*}
&&\varphi_{[A,(0,0)]} = \varphi_{[C_1\cdots C_q,(0,0)]}\\
&=&\varphi_{[C_1,(0,0)]}|[C_2\cdots
C_q,(0,0)]+\varphi_{[C_2,(0,0)]}|[C_3\cdots C_q,(0,0)] + \cdots +
\varphi_{[C_q,(0,0)]},
\end{eqnarray*}
with $q\leq m_1log\mu(A)+m_2$ terms on the right-hand side. By
Lemma \ref{secondlemma}, we have
$$|\varphi_{[A,(0,0)]}(\tau,z)|\leq e^{2\pi
m\frac{y^2}v}K^*\mu(A)^e(|\tau|^\eta+v^{-\eta})q\leq e^{2\pi
m\frac{y^2}v}K_1^*\mu(A)^{e+1}(|\tau|^\eta+v^{-\eta}),$$ where
$\eta = 6e+2k$ and we have used $q\leq m_1log\mu(A)+m_2\leq
m_3\mu(A)$. Lemma \ref{thirdlemma} yields
$$|\varphi_{[A,(0,0)]}(\tau,z)|\leq
e^{2\pi
m\frac{y^2}v}K_2^*(c^2+d^2)^{e+1}(|\tau|^\eta+v^{-\eta}),$$ and,
by Lemma \ref{firstlemma},
\begin{equation}\label{inequality}
|\varphi_{[A,(0,0)]}(\tau,z)|\leq e^{2\pi
m\frac{y^2}v}K_2^*|c\tau+d|^{2e+2}(\frac{1+4|\tau|^2}{v^2})^{e+1}(|\tau|^\eta+v^{-\eta}).
\end{equation}
Note that
$$|\varphi_{[A,(\lambda,0)]}(\tau,z)| \leq
|\varphi_{[A,(0,0)]}|[I,(\lambda,0)](\tau,z)| +
|\varphi_{[I,(\lambda,0)]}(\tau,z)|.$$ Hence, by
(\ref{inequality}) and Lemma \ref{elliptic}
\begin{eqnarray*}
|\varphi_{[A,(\lambda,0)]}(\tau,z)| &\leq& e^{2\pi
m\frac{y^2}v}K_2^*|c\tau+d|^{2e+2}(\frac{1+4|\tau|^2}{v^2})^{e+1}(|\tau|^\eta+v^{-\eta})\\
&\ &+|\lambda|K(|\tau|^\rho+v^{-\sigma})e^{2\pi m\frac{y^2}v}.
\end{eqnarray*}
 To prove the convergence of the series $\Phi(\tau,z)$ we need to
estimate the absolute value of the general term of the series.
This is
\begin{eqnarray*}
&\ &|\varphi_\gamma(\tau,z)(1|_{\omega,k,m}[A,(\lambda,0)])(\tau,z)|\\
&<& e^{2\pi
m\frac{y^2}v}K_2^*|c\tau+d|^{2e+2}(\frac{1+4|\tau|^2}{v^2})^{e+1}(|\tau|^\eta+v^{-\eta})|
(1|_{\omega,k,m}[A,(\lambda,0)])(\tau,z)|\\
&\ & +|\lambda|K(|\tau|^\rho+v^{-\sigma})e^{2\pi
m\frac{y^2}v}|(1|_{\omega,k,m}[A,(\lambda,0)])(\tau,z)|\\
&=& e^{2\pi
m\frac{y^2}v}K_2^*(\frac{1+4|\tau|^2}{v^2})^{e+1}(|\tau|^\eta+v^{-\eta})|(1|_{\omega,k-2e-2,m}[A,(\lambda,0)])(\tau,z)|\\
&\ &+|\lambda|K(|\tau|^\rho+v^{-\sigma})e^{2\pi
m\frac{y^2}v}|(1|_{\omega,k,m}[A,(\lambda,0)])(\tau,z)|,
\end{eqnarray*}
where $A = \sm a&b\\c&d\esm$. By Lemma \ref{convergence} we know that
the series $\Phi(\tau,z)$ converges if $k> 2e+5$ and $k>4$. {\qed}



\begin{thebibliography}{99}




\bibitem{BO1} K. Bringmann and K. Ono,  The $f(q)$ mock theta function
conjecture and partition ranks, Invent. Math. 165 (2006), no. 2,
243--266.


\bibitem{BZ}  K. Bringmann and S.Zwegers, Rank-crank type PDE's and non-holomorphic Jacobi
forms, to appear(2009).


 \bibitem{BO2} J. Bruinier and K. Ono, Heegner divisors, L-functions, and Maass forms,
  to appear in Annals of Mathematics(2009).

\bibitem{C-J} Y. Choie, A short note on the full Jacobi group,
Proceedings of the AMS, Vol 123, No 9, Spe 1995, 2625-2628.

\bibitem{C1} Y. Choie, Half integral weight Jacobi forms and
periods of modular forms, Manuscripta, 104, 124-133(2001).

\bibitem{C2} Y. Choie, Correspondence among Eisenstein series
$E_{2,1}(\tau,z), H_{\frac{3}{2}}(\tau)$ and $E_2(\tau)$,
Manuscripta Math., 93, 177-187 (1997).

\bibitem{C3} Y. Choie, Higher modular symbols and Hecke Operators,
in preparation (2009).


\bibitem{CI-2} Y. Choie and S. Lim, Heat operators, Lerch Sums,
Appell functions and Eichler Integral, Preprint(2009).

\bibitem{CK}Y. Choie and H. Kim, An analogy of Bol's result on
Jacobi forms and Siegel modular forms, Jour of Math. Analysis and
Applications, 257, 79-88(2001).




\bibitem{EZ} M. Eichler and D. Zagier, The Theory of Jacobi forms,
Progress in Mathematics, 55. Birkhuser Boston, Inc., Boston, MA,
1985.


\bibitem{HM} J. Hilgert and D. Mayer, Transfer operators and dynamical zeta
functions for a class of lattice spin models, Comm. Math. Phys.
232 (2002), no. 1, 19--58.


\bibitem{HZ} F. Hirzebruch and D. Zagier,  Intersection Numbers of
curves and Hilbert modular surfaces and Modular forms of
Nebentypus,  Invent. Math, 36, 57-113(1976).


\bibitem{KZ}W. Kohnen and D. Zagier,  Modular forms with rational periods,
Modular forms (Durham, 1983), 197--249, Ellis Horwood Ser. Math. Appl.: Statist.
 Oper. Res., Horwood, Chichester, 1984.


\bibitem{K} M. Knopp, Some New Results on the Eichler
Cohomology of Automorphic Forms, Bull. Amer. Math. Soc.
80(1974),607-632.


\bibitem{K1} M. Knopp,  Rademacher on $J(\tau),$ Poincare series
of nonpositive weights and the Eichler cohomology, Notices Amer.
Math. Soc. 37 (1990), no. 4, 385--393.


\bibitem{K2}M. Knopp, Recent developments in
the theory of rational period
functions,
Number theory (New York, 1985/1988), 111--122,
 Lecture Notes in Math., 1383, Springer, Berlin, 1989.

\bibitem{LZ} J. Lewis and D. Zagier,  Period functions for Maass wave forms. I.
Ann. of Math. (2) 153 (2001), no. 1, 191--258.

\bibitem{M} Y. Manin, Remarks on modular symbols for Maass wave
forms, Arxiv:0803.3270v1 (2008).

\bibitem{M1} Y. Manin and M. Marcolli,   Continued fractions, modular
symbols, and noncommutative geometry. Selecta Math. (N.S.) 8
(2002), no. 3, 475--521.



\bibitem{Mor} L. J. Mordell, The value of the definite integral
$\int_{-\infty}^{\infty} \frac{e^{at^2+bt}}{e^{ct}+d} dt$,
Quarterly J. of Math 68, 1920, 329-342.


\bibitem{Mu} T. Muhlenbruch,   Hecke operators on period functions
for the full modular group. Int. Math. Res. Not. 2004, no. 77, 4127--4145.


\bibitem{P} W. Pribitkin, {\bf{CHECK}}.

\bibitem{P} W. Pribitkin, Eisenstein series and Eichler Integrals,
Contemporary Mathematics, Vol 251,463-467,  2006.

\bibitem{R} S. Ramanujan, The lost notebook and other unpublished
papers, Narosa Publishing House, New Delhi, 1987.





\bibitem{STY} A. M. Semikhatov, A.Taormina and I. Yu. Tipunin,
Higher-Level Appell functions, Modular transformations and
Characters, Comm. Math. Phys. 255 (2005), no. 2, 469--512.

\bibitem{S} A. M. Semikhatov,   Higher string functions, higher-level Appell
functions, and the logarithmic $\widehat{\rm sl}(2)\sb k/{\rm
u}(1)$ CFT model. Comm. Math. Phys. 286 (2009), no. 2, 559--592.



 \bibitem{Za} D. Zagier, Ramanujan's Mock Theta functions and
 their applications, S\'eminaire Bourbaki, 60\'eme ann\'ee, $ 2
006-2007, N^o 986$.

 \bibitem{Ze} S. Zwegers, Mock Theta Functions, PH.D Thesis,
 Universiteit Utrecht, 2002.

 \bibitem{Ze1} S. Zwegers, Mock $\theta$-functions and real
 analytic modular forms, In "$q$-series with Applications to
 Combinatorics, Number Theory and Physics," Contemp. Math. 291,
 Amer. Math. Soc., 2001, 269-277.



\end{thebibliography}
\end{document}